# COMBINATORICS OF RANDOM TENSOR MODELS

Adrian TANASA[1,2]

[1]LIPN, Institut Galilée, CNRS UMR 7030, Université Paris-Nord
99 av. Clement, 93430 Villetaneuse, France
[2]"Horia Hulubei" National Institute for Physics and Nuclear Engineering,
P. O. B. MG-6, 077125 Magurele, Romania
E-mail: adrian.tanasa@ens-lyon.org

In this short overview we introduce group field theory, a particular class of random tensor models, which represents nowadays one of the candidates for a fundamental theory of quantum gravity. We insist on the combinatorial richness of associated structures, namely tensor graphs, natural generalization of ribbon graphs (or combinatorial maps).

*Key words*: Feynman tensor graphs, bubbles, group field theory, (random) tensor models.

Random tensor models can be seen as a generalization of random matrix models and they are nowadays one of the candidates for a fundamental theory of quantum gravity. The theoretical physics framework is a quantum field theoretical (QFT) one, namely group field theory (GFT) (the interested reader can turn to the review papers [1] or [2]). Roughly speaking, the main idea is that if matrix models can be seen as dual to triangulations of two-dimensional surfaces, this can then be extended – using tensor models instead of matrix models – to three-dimensional and finally to four-dimensional spaces. It is worth emphasizing, that unlike string theories, these models do not require extra-dimensions, so that one can work in a space with dimension four, the dimension of space-time at our energy scale.

Let us also mention that an important regain of interest for GFTs is to be noticed lately in the mathematical physics community.

In the combinatorial simplest case, the elementary cells that, by gluing together form the space itself, are the $D$-simplices ($D$ being the dimension of space, here taken equal to three and four, as already mentioned above). Since a $D$-simplex has $(D+1)$ facets on its boundary, the backbone of group field theoretical models in $D$-dimension should be some abstract $\Phi^{D+1}$ interaction on rank $D$ tensor fields $\Phi$ (see Fig. 1, Fig. 2, and respectively Fig. 3 which represent, for $D = 2$, $D = 3$, and respectively $D = 4$ vertices of these Feynman quantum gravity tensor graphs).

GFT is defined as a quantum field theory on group manifolds. Its associated Feynman graphs are, for the three-dimensional case, rank three tensor graphs and for the four-dimensional case, rank four tensor graphs.

In order to define such a scalar Feynman graph one needs a propagator (associated to the edge of the graph) and an interaction (associated to the graph vertex). The edge is graphically given by Fig. 4. Each such GFT tensor graph edge has, in three dimensions, three strands and respectively four strands in four dimensions.

As already mentioned above, GFT graphs are dual to triangularizations of space-time. Thus, in the simplest two-dimensional case, the building block of such a triangularization is of course a triangle. One can now take a dual point of view and encode the information coming from this triangularization in a ribbon graph picture (see for example the two-dimensional triangularization of Fig.1).

The appropriate valence of the associate ribbon graph vertex is thus three: the vertex (Fig.1) corresponds to the triangle (the building block of the triangularization) and the three incoming/outgoing edges correspond to the three edges of the original triangle.

The two strands of each edge correspond to the two vertices joined by the respective edge in the initial triangle. For the sake of completeness, let us also mention that this point of view relates to the one of two-dimensional quantum gravity matrix models.

The situation is similar for the three-dimensional case. The valence four vertex corresponds to a tetrahedron, the simplest building-block of a three-dimensional space-time. The four triangles constructing the tetrahedron correspond to the four edges intersecting at the respective vertex (see for example Fig. 2).

The three strands of such an edge correspond to the three edges of the respective triangle (face of the tetrahedron), thus generalizing the image of the two-dimensional case.

For the sake of simplicity, we analyze in this paper orientable three-dimensional tensor graphs. This means that the simplicial complexes dual to these graphs are orientable (the interested reader can turn for example to [3] for more details on this point). Let us also mention that the concepts used in this paper do not seem to apply directly to general tensor graphs.

The strand structure described above makes these theories extremely rich from a combinatorial and topological point of view. Thus, in addition to faces (closed circuits) which are natural to define for ribbon graphs, one can further generalize this topological notion for the case of tensor graphs.

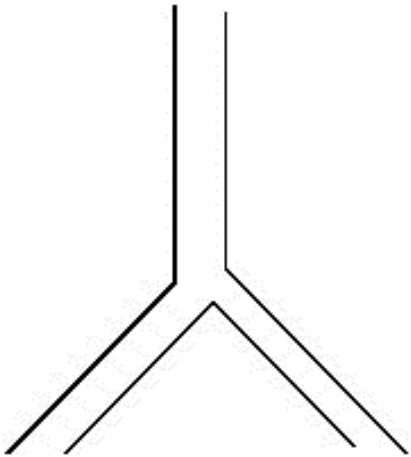
Fig. 1 – The valence three vertices of two-dimensional GFT.

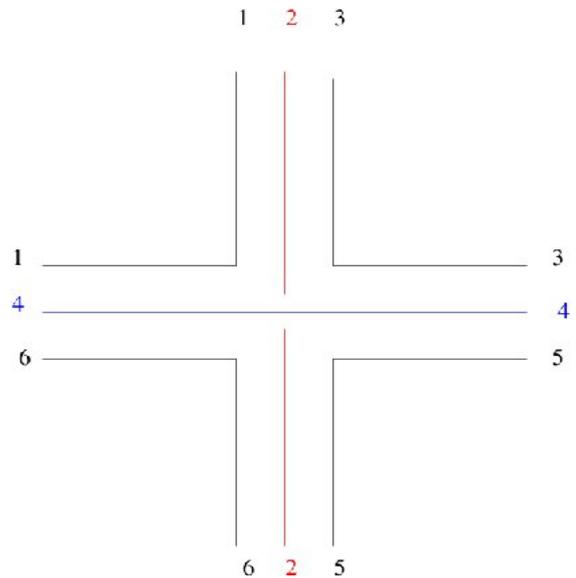
Fig. 2 – The valence four vertices of three-dimensional GFT.

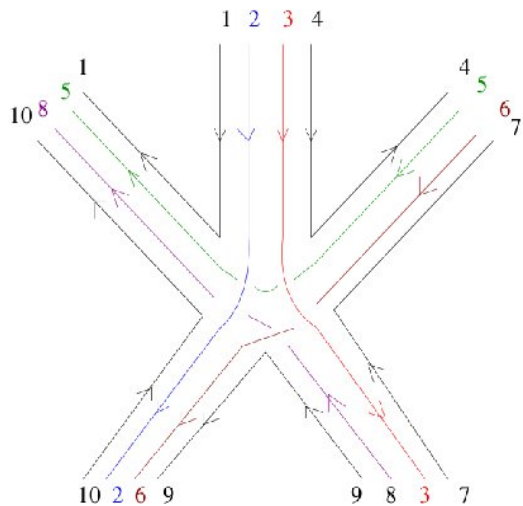
Fig. 3 – A vertex of a four-dimensional GFT. We have chosen here a particular matching and orientation for each of the strands.

Such a generalization is the notion of *bubble*, which can be defined as a closed three-dimensional region of the graph.

A bubble can thus be represented as a ribbon graph, and thus, from a topological point of view, it is natural to define its Euler characteristic $\chi$ (and hence its genus). This is just a translation of the fact that such a ribbon graph can be drawn on a two-dimensional manifold (each graph defining the surface on which it is drawn).

One thus speaks of *planar* or *non-planar bubbles*. Let us give some more explanations on this. For the example of Fig. 5, the bubbles are given in Fig. 6 and Fig. 7 (note that a slightly different convention for drawing the vertices has been taken; nevertheless this does not interfere with the presentation of the combinatorial notions of this paper). The first of them is non-planar, while the second of them is planar.

Note that in [3] a precise algorithm for characterizing these bubbles was given for an arbitrary GFT tensor graph. A cohomological definition of this notion of bubble was given in [4].

As faces do for ribbon graphs, bubbles are known to play a crucial rôle for these tensor graphs. Thus, if in QFT the divergences (which require renormalization) are related to *cycles* (known under the name of *loops* in mathematical physics), in GFT these divergences are related to this new notion of bubbles.

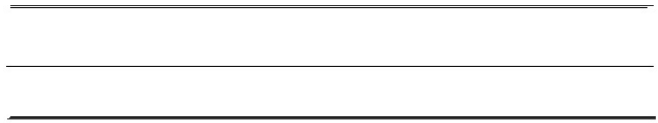

Fig. 4 – The edge of the GFT tensor graphs.

Lately, a particular class of tensor graphs – the *colorable tensor graphs*, has been intensively studied (see the review paper [5] and references within).

In three dimensions, a **colorable graph** is a graph where each edge can be assigned a certain color belonging to a finite set of four elements, such that each of the four incoming/outgoing edges of an arbitrary vertex has a distinct color and at each vertex one has a clockwise or an anti-clockwise cyclic ordering for the colors of the four incoming/outgoing edges.

As a direct consequence of this definition, one can conclude that a colorable graph has two types of vertices: a white (or positive)} one (where the order is anti-clockwise) one and a black (or negative) one (where the order is clockwise). When a graph has two types of vertices, it is called in graph theory a *bipartite graph*. Furthermore, any edge of such a graph connects a positive to a negative vertex.

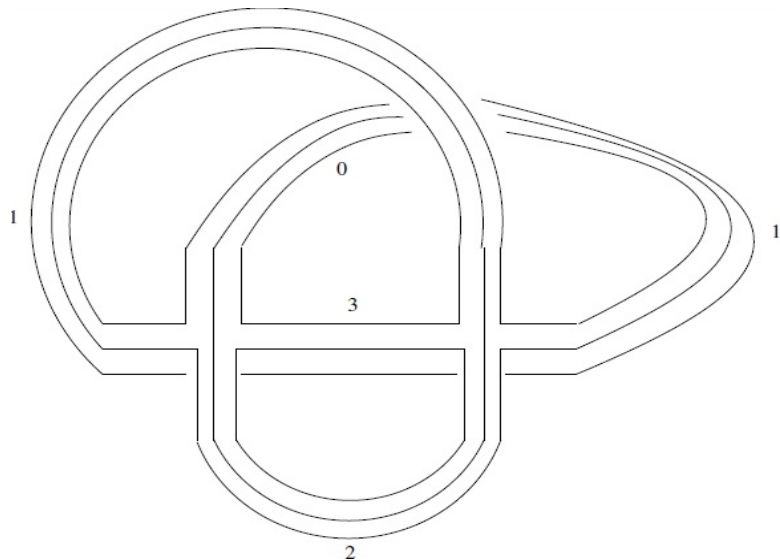

Fig. 5 – An example of a GFT graph. We denote by $e_0, \ldots, e_3$ the four edges.

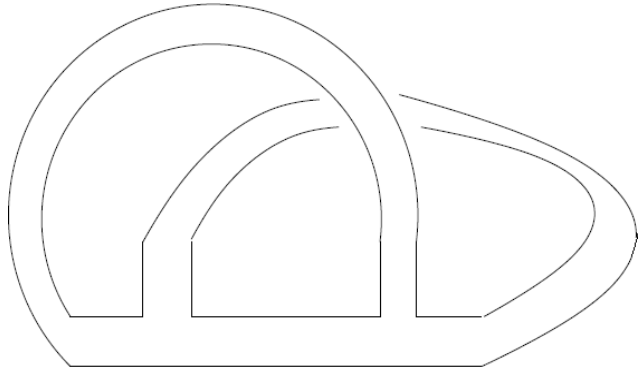

Fig. 6 – One of the bubbles of the graph of Fig. 5. It has two vertices,
three edges and one face and thus genus 1; it is non-planar.

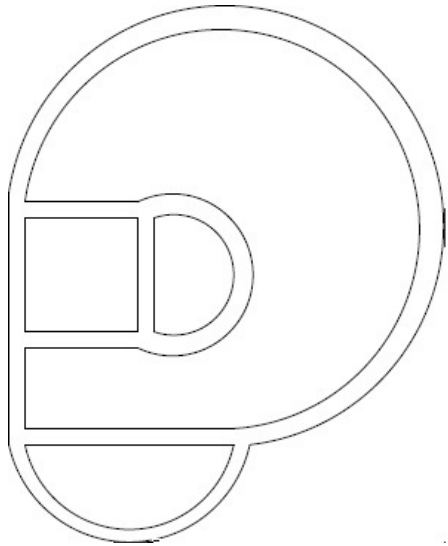

Fig. 7 – The second of the bubles of the graph of the Fig. 5;
it has genus 0, being hence planar.

Furthermore, let us emphasize on the following. When drawing these colorable graphs, one can drop the cumbersome stranded standard picture, because all the information is coded in the coloring. This can actually be done for any GFT tensor graph (colorable or not) which does not allow twists in the edges (or vertices). In QFT language the coloring presented here is related to the fact that one deals with a complex field.

As expected, the topological notion of bubbles is defined in an easier way for these colorable graphs. Thus, a bubble can now be defined as a connected component of a subgraph formed only of edges of colors belonging to a subset of cardinal three of the initial set of colors. Let us end this brief overview by mentioning that alternative models to colorable tensor models have been recently proposed, the **multi-orientable models** [6].

The main idea is that one imposes not only that there are no twists between the various strands of an edge, but also that a regular vertex has two positive and two negative corners (following in a cyclic order) and that an edge can connect only a positive to a negative corner. This idea, which was proved successful in noncommutative quantum field theory, leads to the fact that a smaller class of graphs is discarded, when comparing the situation with the colorable models one. From a QFT point of view, the multi-orientable models also need a complex scalar field.

It is thus an appealing perspective for future work to investigate weather or not the various theoretical physics developments achieved for colorable models (see again [5]), can be adapted for these alternative models.

## ACKNOWLEDGEMENTS

The CNRS PEPS "CombGraph", PN 09370102, CNCSIS Tinere echipe 77/04.08.2010 and Idei 454/2009 are acknowledged.